\documentclass{amsart}

\usepackage{amsfonts, amsmath, amssymb, amsthm, amsxtra, latexsym}
\usepackage[all]{xy}

\theoremstyle{plain}

\newtheorem{theorem}{Theorem}[section]
\newtheorem{proposition}[theorem]{Proposition}
\newtheorem{lemma}[theorem]{Lemma}
\newtheorem{corollary}[theorem]{Corollary}

\theoremstyle{definition}

\newtheorem{definition}[theorem]{Definition}

\newtheorem{remark}[theorem]{Remark}

\numberwithin{equation}{section}

\newcommand{\rep}{{\pi}}
\newcommand{\hilb}{{H}}
\newcommand{\bounded}{{B(\hilb)}}
\newcommand{\borelfunc}{{\mathcal{B}_{b}(X)}}

\newcommand{\res}{{\textup{Res}}}

\begin{document}


\title[Commutative ${C^*}$-algebras and {\Large $\sigma$}-normal morphisms]{Commutative $\mathbf{C^*}$-algebras and {\LARGE $\mathbf{\sigma}$}-normal morphisms}


\translator{}

\dedicatory{}

\author[Marcel~de~Jeu]{Marcel~de~Jeu}

\begin{abstract}
We prove in an elementary fashion that the image of a commutative monotone
$\sigma$-complete $C^*$-algebra under a $\sigma$-normal morphism is again
monotone $\sigma$-complete and give an application of this result in spectral theory.

\end{abstract}

\date{}



\thanks{During the preparation of this paper the author was supported by a PIONIER grant of the Netherlands Organisation for Scientific Research (NWO)}

\address{M.F.E.~de~Jeu\\
         Mathematical Institute\\
         Leiden University\\
         P.O. Box 9512\\
     2300 RA Leiden\\
     The Netherlands}



\email{mdejeu@math.leidenuniv.nl}

\urladdr{}

\maketitle

\section{Introduction}

The image of a monotone $\sigma$-complete $C^*$-algebra under a
$\sigma$-normal morphism is again monotone $\sigma$-complete, as can be deduced from somewhat technical results in the literature \cite{Christensen, Wright}. In Section~\ref{sec:mainresult} of the present paper, however, we give  a straightforward proof when the algebra is
commutative, using only the commutative Gelfand--Naimark theorem. We also indicate how the general statement can be inferred. Section~\ref{sec:application} contains an application in
spectral theory which justifies special interest in the commutative case.

\section{Main result}\label{sec:mainresult}

The $C^*$-algebras in
this paper are not necessarily unital and neither are (if applicable) the
morphisms. This being said, we recall the relevant order definitions.
\begin{definition}
Cf.\ \cite[3.9.2]{Pedersen}.
\begin{enumerate}
\item A $C^*$-algebra $A$ is \emph{monotone $\sigma$-complete} if every bounded increasing sequence of self-adjoint elements of $A$ has a supremum in $A$. \item A $C^*$-subalgebra $A$ of a
$C^*$-algebra $B$ is a \emph{monotone $\sigma$-closed $C^*$-subalgebra of $B$} if $\sup^B_{n\geq 1} a_n\in A$, whenever $a_1\leq a_2\leq\ldots$ is a bounded increasing sequence of self-adjoint
elements of $A$ which has a supremum $\sup^B_{n\geq 1} a_n$ in $B$. \item A morphism $\phi: A\mapsto B$ between two $C^*$-algebras is a \emph{$\sigma$-normal morphism} if $\phi(\sup_{n\geq
1}a_n)=\sup_{n\geq 1}\phi(a_n)$, for each bounded increasing sequence $a_1\leq a_2\leq\ldots$ of self-adjoint elements of $A$ which has a supremum in $A$.
\end{enumerate}
\end{definition}

Our main result then reads as follows.

\begin{theorem}\label{thm:maintheorem}
Let $A$ be a commutative monotone $\sigma$-complete $C^*$-algebra. Suppose $\phi:A\mapsto B$ is a $\sigma$-normal morphism into a $C^*$-algebra $B$. Then $\phi(A)$ is a commutative monotone $\sigma$-complete $C^*$-algebra, and a monotone $\sigma$-closed $C^*$-subalgebra of $B$.
\end{theorem}

The proof of Theorem~\ref{thm:maintheorem} which we will now give uses no
result beyond the commutative Gelfand--Naimark theorem. We start with an
isometric $\sigma$-monotone lifting result for functions.

\begin{lemma}\label{lem:liftingfunctions} Let $X$ be a topological space, and let $C_0(X)$ denote the continuous functions on $X$ vanishing at infinity. Suppose that $Y\subset X$ is a nonempty subset.
\begin{enumerate}
\item If $0\leq f_1\leq f_2\leq\ldots$ is a sequence of functions on $Y$, such that each $f_n$ is the restriction of some element of $C_0(X)$ to $Y$, then there exists a sequence $0\leq g_1\leq
g_2\leq\ldots$ in $C_0(X)$ such that, for each $n$, $g_n$ restricts to $f_n$ and $\Vert g_n\Vert_\infty=\Vert f_n\Vert_\infty$.
\item If $f_1\geq f_2\geq\ldots\geq 0$ is a sequence of functions on
$Y$, such that each $f_n$ is the restriction of some element of $C_0(X)$ to $Y$, then there exists a sequence $g_1\geq g_2\geq\ldots\geq 0$ in $C_0(X)$ such that, for each $n$, $g_n$ restricts to
$f_n$ and $\Vert g_n\Vert_\infty=\Vert f_n\Vert_\infty$.
\end{enumerate}
\end{lemma}

\begin{proof}
Suppose that $h_n\in C_0(X)$ restricts to $f_n$. Replacing $h_n$ with $|h_n|$, we may assume that $h_n\geq 0$. After a  subsequent replacement of $h_n$ with $\min (h_n,\Vert f_n\Vert_\infty)\in
C_0(X)$, we may assume that $h_n\geq 0$ and that $\Vert h_n\Vert_\infty=\Vert f_n\Vert_\infty$. In the case of an increasing sequence, define $g_n=\max_{1\leq i\leq n} h_i$. In the case of a
decreasing sequence, define $g_n=\min_{1\leq i\leq n} h_i$. Then the $g_n$ have the required properties.
\end{proof}

Conbining this with the commutative Gelfand--Naimark theorem, we obtain the
following generalization of the above result. It is somewhat stronger than
needed for our purpose, as we will use only the first part without the isometric property in the
sequel.

\begin{proposition}\label{prop:liftingelements}
Let $A$ and $B$ be commutative $C^*$-algebras. Suppose that $\phi:A\mapsto B$ is a surjective morphism.
\begin{enumerate}
\item If $0\leq b_1\leq b_2\leq\ldots$ is a sequence in $B$, then there exists a sequence $0\leq a_1\leq a_2\leq\ldots$ in $A$, such that, for each $n$, $\phi(a_n)=b_n$ and $\Vert a_n\Vert=\Vert
b_n\Vert$. \item If $b_1\geq b_2\geq\ldots\geq 0$ is a sequence in $B$, then there exists a sequence $a_1\geq a_2\geq\ldots\geq 0$ in $A$, such that, for each $n$, $\phi(a_n)=b_n$ and $\Vert
a_n\Vert=\Vert b_n\Vert$.
\end{enumerate}
\end{proposition}

\begin{proof}
We may assume that $B\neq 0$. In that case, define the canonical map
$\phi^*:\widehat B\mapsto\widehat A$ between the spectra by
$\phi^*(\beta)=\beta\circ\phi$ for $\beta\in \widehat B$. Since $\phi$ is
surjective, $\phi^*$ is injective. We identify $\widehat B$ as a set with its
image in $\widehat A$ using $\phi^*$, and we thus view $C_0(\widehat B)$ as
functions on $\widehat B\subset\widehat A$. Actually, $\phi^*$ is a
homeomorphic embedding of $\widehat B$ into $\widehat A$, but this is not
needed. Using these two identifications, the following diagram is then
commutative:
\[
\xymatrix{A \ar[r]^\phi \ar[d]_{\Gamma_A} &B \ar[d]^{\Gamma_B}\\
C_0(\widehat A) \ar[r]_{\res}&C_0(\widehat B)}
\]
Here $\Gamma_A$ and $\Gamma_B$ denote the Gelfand--Naimark isomorphisms, and $\res$ is the restriction mapping. The proposition now follows from an application of Lemma~\ref{lem:liftingfunctions} at
the bottom part of the diagram.
\end{proof}

\begin{proposition}\label{prop:finalprop}
Let $A$ be a commutative monotone $\sigma$-complete $C^*$-algebra. Suppose $\phi:A \mapsto B$ is a $\sigma$-normal morphism into a $C^*$-algebra $B$. Let $b_1\leq b_2\leq\ldots$ be a bounded sequence of self-adjoint elements of $\phi(A)$. Then this sequence has a supremum in $B$, and this supremum in $B$ is an element of $\phi(A)$.
\end{proposition}

\begin{proof}
We may assume that $b_1\geq 0$. Using Proposition~\ref{prop:liftingelements} for the morphism $\phi: A\mapsto\phi(A)$, we find a bounded sequence $a_1\leq a_2\leq\ldots$ of self-adjoint elements of
$A$, such that $\phi(a_n)=b_n$ for all $n$. Since $A$ is monotone $\sigma$-complete, $\sup_{n\geq 1}a_n$ exists in $A$. By the $\sigma$-normality of $\phi$, $\sup_{n\geq 1} b_n$ exists in $B$,
since it is equal to $\phi( \sup_{n\geq 1}a_n)$. This also shows that this supremum in $B$ is an element of $\phi(A)$.
\end{proof}

Proposition~\ref{prop:finalprop} implies Theorem~\ref{thm:maintheorem}, as the reader will easily verify.

\begin{remark}\label{rem:comparison}
As mentioned in the Introduction, Theorem~\ref{thm:maintheorem} also holds
in the noncommutative case. In fact, the noncommutative version of
Proposition~\ref{prop:finalprop}
--- which implies the noncommutative version of Theorem~\ref{thm:maintheorem}
--- is valid, as can be seen from \cite[Lemma~1.1]{Wright}, which is in turn
based on \cite[Proposition~5]{PedersenDecomposition}. The proof of the latter
basic result relies on the general Gelfand--Naimark theorem and
involves some well chosen use of the Borel functional calculus.

As an alternative route, in the unital noncommutative case the monotone
$\sigma$-completeness of $\phi(A)$ in Theorem~\ref{thm:maintheorem} also
follows from \cite[Proposition~1.3]{Wright}, which is in turn based on
\cite[Lemma~2.13]{Christensen}. A closer inspection of the proof of the latter
result shows that it in fact yields Proposition~\ref{prop:finalprop} in the
unital noncommutative setting, and therefore Theorem~\ref{thm:maintheorem} in
the unital noncommutative case again. This key proof of
\cite[Lemma~2.13]{Christensen} relies itself on the fact that a unital
monotone $\sigma$-complete $C^*$-algebra is a Rickart $C^*$-algebra in which
the projections form a lattice. As remarked in the introduction of \cite{Christensen}, this follows from \cite{KadisonPedersen}.

The above makes it clear that in the commutative case the proof of
Theorem~\ref{thm:maintheorem} as given in this section is considerably more elementary than the alternative approach of establishing and subsequently specializing the noncommutative version.

\end{remark}

\section{Application in spectral theory}\label{sec:application}

To conclude, we give an application of Theorem~\ref{thm:maintheorem} in
spectral theory that justifies special interest in the commutative case. It is
based on the following.

\begin{corollary}\label{cor:maincor}Let $A$ be a commutative monotone $\sigma$-complete $C^*$-algebra and suppose
$\rep:A\mapsto\bounded$ is a $\sigma$-normal representation in a separable
Hilbert space $\hilb$. Then $\rep(A)$ is strongly closed. If $A$ and $\rep$ are
unital, then $\rep(A)$ is a von Neumann algebra.
\end{corollary}

Indeed, $\rep(A)$ is a monotone $\sigma$-closed $C^*$-subalgebra of $\bounded$
by Theorem~\ref{thm:maintheorem}. Hence the Up--Down theorem for separable
Hilbert spaces \cite[Theorem 2.4.3]{Pedersen} implies that it is strongly
closed, establishing the corollary (which is also valid in the noncommutative case in view of Remark~\ref{rem:comparison}).

As an application of Corollary~\ref{cor:maincor} in spectral theory, let $X$ be
a compact Hausdorff space with associated $C^*$-algebras $C(X)$ of continuous
functions and $\borelfunc$ of bounded Borel measurable functions. Suppose that
$\rep: C(X)\mapsto\bounded$ is a unital representation. Then $\rep$ extends uniquely to a unital representation of
$\borelfunc$, again denoted by $\rep$, such that
\begin{equation}\label{eq:borelrep}
(\rep (f)\xi,\xi)=\int_X f(x)\, d\mu_\xi(x)\quad(f\in\borelfunc,\,\xi\in\hilb).
\end{equation}
Here $\mu_\xi$ denotes the positive and bounded unique regular Borel measure on
$X$ which is provided by the Riesz representation theorem when one requires
\eqref{eq:borelrep} to hold for all $f\in C(X)$.

It is easy to see that $\rep$ maps $\borelfunc$ into the strong closure of $\rep(C(X))$, 
implying that $\rep(\borelfunc)$ and  $\rep(C(X))$ have the same strong
closure. If $\hilb$ is separable, then $\rep(\borelfunc)$ is actually equal to
the strong closure of $\rep(C(X))$ \cite[proof of Proposition~9.5.3]{KadRing}. This
description of the von Neumann algebra generated by $\rep(C(X))$ as being equal
to $\rep(\borelfunc)$ is a basic ingredient for further investigation of
the unital separable representation $\rep$ of $C(X)$.

The proof of the equality of $\rep(\borelfunc)$ and the strong closure of
$\rep(C(X))$ in the separable case is usually based on measure-theoretical
arguments and is somewhat more involved than others in this circle of
ideas, cf.~\cite[proof of Proposition~9.5.3]{KadRing}, but Corollary~\ref{cor:maincor}, which is based on order properties, provides an
alternative approach. Indeed, $\borelfunc$ is monotone $\sigma$-complete and,
as a consequence of \eqref{eq:borelrep} and the monotone convergence theorem,
$\rep$ is a $\sigma$-normal representation of $\borelfunc$. Therefore, by
Corollary~\ref{cor:maincor}, $\rep(\borelfunc)$ is strongly closed, hence equal
to the strong closure of $\rep(C(X))$.

In fact, by the same reasoning Corollary~\ref{cor:maincor} also implies ---
still for a unital separable representation --- that the image of the monotone
$\sigma$-completion of $C(X)$ in $\borelfunc$ is strongly closed, hence already
equal to the von Neumann algebra generated by $\rep(C(X))$. If $X$ is second
countable then this monotone $\sigma$-completion of $C(X)$ in $\borelfunc$ ---
the Baire algebra --- coincides with $\borelfunc$, but in other cases it may be
a proper $C^*$-subalgebra of $\borelfunc$ \cite[Remarks~6.2.10]{PedersenNow}.





\begin{thebibliography}{99}

\bibitem{Christensen} E.~Christensen, \emph{Non commutative integration for monotone sequentially closed $C^*$-algebras}, Math.\ Scand.\ {\bf 31} (1972), 171-190.

\bibitem{KadisonPedersen} R.V.~Kadison and G.K.~Pedersen, \emph{Equivalence in operator algebras}, Math.\ Scand.\ {\bf 27} (1970), 205-222.

\bibitem{KadRing} R.V.\ Kadison and J.R.\ Ringrose, Fundamentals of the Theory of Operator Algebras, Volume II, Academic Press, London, 1986.

\bibitem{PedersenDecomposition} G.K.~Pedersen, \emph{A decomposition theorem for $C^*$-algebras}, Math.\ Scand.\ {\bf 22} (1968), 266-268.

\bibitem{PedersenNow} G.K.\ Pedersen, Analysis Now, Springer-Verlag, New York, 1989.

\bibitem{Pedersen} G.K.\ Pedersen, $C^*$-Algebras and their automorphism groups, Academic Press, London, 1979.

\bibitem{Wright} J.D.M.~Wright, \emph{On minimal $\sigma$-completions of $C^*$-algebras}, Bull.\ London Math.\ Soc.\ {\bf 6} (1974), 168-174.

\end{thebibliography}
\end{document}